\documentclass{article}
\usepackage{graphicx}
\usepackage{psfrag}
\usepackage{psfig}

\newcommand{\half}{{\textstyle{\frac{1}{2}}}}
\newcommand{\onethird}{{\textstyle{\frac{1}{3}}}}

\begin{document}   

\begin{center}
\begin{LARGE}
\textbf{A Fast Solver for Systems of Reaction-Diffusion Equations} \\[2ex]
\end{LARGE}
\begin{Large}
M.~Garbey,\footnote{
     CDCSP-ISTIL, Universit\'{e} Claude Bernard Lyon 1,
     cedex 69622 Villeurbanne, France}
H.~G.~Kaper,\footnote{
     Mathematics and Computer Science Division,
     Argonne National Laboratory, Argonne, IL 60439, USA.
     Work supported by the Mathematical, Information, and
     Computational Sciences Division subprogram of the Office of Advanced
     Scientific Computing Research, U.S.~Department of Energy, under
     Contract W-31-109-Eng-38.}
and
N.~Romanyukha\footnote{
     Institute for Mathematical Modeling,
     Russian Academy of Science, 125047, Moscow, Russia.
     Supported by the Russian Foundation of Basic Research
     under Grant 00-01-0291.}
\end{Large}
\end{center}

\medskip

\section{Introduction}
In this paper we present a fast algorithm for
the numerical solution of systems of
reaction-diffusion equations,
\begin{equation} \label{GKR-global}
  \partial_t u
  + a \cdot \nabla u
  = \Delta u
  + F (x, t, u) , \quad
  x \in \Omega \subset \mathbf{R}^3, \, t > 0 .
\end{equation}
Here, $u$ is a vector-valued function,
$u \equiv u(x, t) \in \mathbf{R}^m$,
$m$ is large,
and the corresponding system of ODEs,
$\partial_t u = F(x, t, u)$, is stiff.
Typical examples arise in air pollution studies,
where $a$ is the given wind field and
the nonlinear function $F$ models the
atmospheric chemistry.

The time integration of Eq.~(\ref{GKR-global})
is best handled by the method of
characteristics~\cite{GKR-Pironneau}.
The problem is thus reduced to designing
for the reaction-diffusion part a fast solver
that has good stability properties for
the given time step and does not require
the computation of the full Jacobi matrix.

An operator-splitting technique,
even a high-order one,
combining a fast nonlinear ODE solver
with an efficient solver for the
diffusion operator is less effective
when the reaction term is stiff.
In fact, the classical Strang splitting method
may underperform a first-order source splitting
method~\cite{GKR-Sportisse}.
The algorithm we propose in this paper
uses an a posteriori filtering technique
to stabilize the computation of the diffusion term.
The algorithm parallelizes well,
because the solution of the large system
of ODEs is done pointwise;
however, the integration of the chemistry
may lead to load-balancing problems~\cite{GKR-Dabdub, GKR-Elbern}.
The Tchebycheff acceleration technique
proposed in~\cite{GKR-Droux} offers
an alternative that complements
the approach presented here.

To facilitate the presentation,
we limit the discusssion to domains $\Omega$
that either admit a regular discretization grid
or decompose into subdomains that admit
regular discretization grids.
We describe the algorithm
for one-dimensional domains in Section~2
and for multidimensional domains in Section~3.
Section~4 briefly outlines future work.

\section{One-Dimensional Domains}
Consider the scalar equation
\begin{equation} \label{GKR-scalar_reaction_diffusion}
  \partial_t u
  =
  \partial_x^2 u
  + f (u) , \quad
  x \in (0, \pi), \,  t > 0 .
\end{equation}
We combine a backward Euler approximation in time
with an explicit finite-difference approximation
of the diffusive term,
\begin{equation} \label{GKR-schema}
  \frac{3 u^{n+1} - 4 u^{n} + u^{n-1}}{2 \Delta t}
  =
  2 D_{xx} u^{n} - D_{x} u^{n-1} + f (u^{n+1}) .
\end{equation}
This scheme is second-order accurate in both space and
time~\cite{GKR-Petzold,GKR-Sandu,GKR-Verwer}.
To analyze its stability, we take the Fourier
transform of the linear equation,
\begin{equation} \label{GKR-Fourier}
  \frac{3 \hat{u}^{n+1} - 4 \hat{u}^{n} + \hat{u}^{n-1}}{2 \Delta t}
  = \Lambda_k (2 \hat{u}^{n} - \hat{u}^{n-1}) ,
\end{equation}
where $\Lambda_k = 2 h^{-2} (\cos(hk) - 1)$,
from which we obtain the stability condition
\begin{equation} \label{GKR-wave}
  2 \frac{\Delta t}{h^2}
  \left| \cos \left( \frac{k \pi}{N} \right)-1 \right|
  < \frac{4}{3} , \quad h = \frac{\pi}{N} .
\end{equation}
Thus we conclude that the time step must satisfy the constraint
\begin{equation}
  \Delta t < \onethird h^2 .
\end{equation}
However, this constraint is imposed by
the high frequencies, which are poorly handled
by second-order finite differences anyway.
For example, with central differences,
the relative error for high-frequency waves
$\cos (kx)$ with $k \approx N$
can grow at a rate of up to $9\%$.
The idea is therefore to relax the constraint
on the time step by applying a filter
after each time step, which removes
the high frequencies but maintains
second-order accuracy in space.    

\subsection{Filters}
By a \emph{filter of order $p$} we mean
an even function
$\sigma : \mathbf{R} \to \mathbf{R}$
that satisfies the conditions
(i)~$\sigma(0)=1$,
(ii)~$\sigma^{(l)}(0)=0$ for $l = 1, \ldots\,, p-1$,
(iii)~$\sigma(\eta)=0$ for $| \eta | \geq 1$, and 
(iv)~$\sigma \in C^{p-1} (\mathbf{R})$.

\smallskip\noindent
\textbf{Theorem}~\cite{GKR-Gottlieb}.
Let $f$ be a piecewise $C^p$ function
with one point of discontinuity, $\xi$,
and let $\sigma$ be a filter of order $p$.
For any point $y \in [0 , 2 \pi]$,
let $d(y) = \min \{| y - \xi + 2 k \pi | : k = -1, 0, 1 \}$.
If
$f^{\sigma}_{N}
  =
  \sum_{k=-\infty}^{\infty}
  {\hat{f}}_k \sigma(k/N) \mathrm{e}^{i k y}$,
then
\[
  \left| f(y) - k/N \right|
  \leq
  C N^{1-p} (d(y))^{1-p} K(f)
  + C N^{1/2 - p} \|f^{(p)}\|_{L^2} ,
\]
where 
\[
  K(f) = \sum_{l=0}^{p-1}
  (d(y))^l |f^{(l)} (\xi^{+}) - f^{(l)}(\xi^{-})|
  \int_{-\infty}^{\infty}
  |G_l^{(p-l)}(\eta)| \, \mathrm{d} \eta .
\]

In other words, a discontinuity of $f$ leads to
a Fourier expansion with an error that is
$O(1)$ near the discontinuity and
$O(N^{-1})$ away from the discontinuity.
We must therefore apply a shift and
extend to $[0, 2\pi]$
before applying a filter.

\subsection{The Algorithm}
We now describe the postprocessing algorithm
that is to be applied after each time step.
(We do not explicitly indicate
the dependence of $u$ on the time step,
and we use the abbreviations
$u_0 = u(0)$ and $u_\pi = u(\pi)$.)

First, we apply a low-frequency shift,
\begin{equation} \label{GKR-step1}
  v(x)
  = u(x)
  - (\alpha_1 + \alpha_2 \cos(x)) , \quad
  \alpha_1 = \half (u_0 - u_{\pi}) , \,
  \alpha_2 = \half (u_0 + u_{\pi}) .
\end{equation}
Then we extend $v$ to $(0,2\pi)$,
using the definition
\begin{equation}
  v(2\pi - x) = - v(x) , \quad x\in(0,\pi) .
\end{equation}
Thus, $v$ is a $2\pi$-periodic function
in $C^1(0,2 \pi)$.
Let $\hat{v}_k$ be the $k$th coefficient
of its Fourier expansion.

Next, we apply an eighth-order filter~\cite{GKR-Gottlieb},
\begin{equation} \label{GKR-step2}
  \sigma_N v(x)
  = \sum_k
  \sigma \left( \kappa \frac{k}{N} \right)
  \hat{v}_k \mathrm{e}^{i k x} ,
\end{equation}
where
\begin{equation}  \label{GKR-filter}
  \sigma(\xi) = (35 - 84 y + 70 y^2 - 20 y^3) y^4 , \quad
  y \equiv y(\xi) = \half (1 + \cos (\pi \xi)) .
\end{equation}
Here, $\kappa$ is a stretching factor, $\kappa > 1$.
The correct choice of $\kappa$ follows
from a Fourier analysis of Eq.~(\ref{GKR-wave}),
\begin{equation} \label{GKR-kappa_limit}
  \kappa > \kappa_c
  = \frac{\pi}{\cos (1 - 2h^2/(3\Delta t))} .
\end{equation}
The choice $\kappa = \frac{1}{2} \kappa_c$
gives satisfactory results,
but in principle one can compute
the optimum value of $\kappa$
at each time step by monitoring
the growth of the high-frequency waves
that have not been completely filtered out.

Finally, we recover $u$ from the inverse shift,
\begin{equation} \label{GKR-step3}
  u(x) = \sigma_N v(x) + \alpha_1 + \alpha_2 \cos (x) .
\end{equation}

The theorem quoted in the preceding section shows that
the filtering process may affect the spatial accuracy
of the method.
Since the filter is applied to a
$2 \pi$-periodic function that is
$C^1$ at the points $x_k = k \pi$, $k \in \mathbf{Z}$,
and $C^2$ everywhere else,
the error is of the order of $N^{-2}$
in the neighborhood of $x_k$
and $N^{-3}$ away from $x_k$.
In principle, we maintain therefore 
second-order accuracy in space
as long as $\kappa$ is of order one.

If the solution
is in $C^3(0,\pi)$ at each time level,
we can improve the algorithm
by replacing the first-order shift~(\ref{GKR-step1})
by a third-order shift,
\begin{equation} \label{GKR-shift4}
  v(x) = u(x) - \sum_{j=1}^4 \alpha_j \cos((j-1) x) ,
\end{equation}
such that the extension of $v$ to a $2 \pi$-periodic
function is in $C^3(0,2 \pi)$.
The first- and third-order derivatives of $v$
are zero at the points $x_k$,
and the second-order derivative
is approximately given by
\begin{equation} \label{GKR-uxx}
  u_{xx}(x_k)
  \approx
  \frac{3 u^{n+1}(x_k) - 4 u^{n}(x_k) + u^{n-1}(x_k)}{2 \Delta t}
  -
  f (u^{n+1}(x_k)) .
\end{equation}
The coefficients $\alpha_j$
are found by solving a linear system of equations,
\begin{eqnarray*}
  \alpha_1 + \alpha_2 + \alpha_3 + \alpha_4 &=& u(0) , \\
  \alpha_1 - \alpha_2 + \alpha_3 - \alpha_4 &=& u(\pi) , \\
           - \alpha_2 - 4 \alpha_3 - 9 \alpha_4 &=& u_{xx}(0) , \\
             \alpha_2 - 4 \alpha_3 + 9 \alpha_4 &=& u_{xx}(\pi) .
\end{eqnarray*}
The third-order shift improves the performance
of the filter for large $\kappa$ and allows
for a larger time step.

\subsection{Numerical Results}
\begin{figure}[htb]
  \begin{center}
  \mbox{\psfig{file=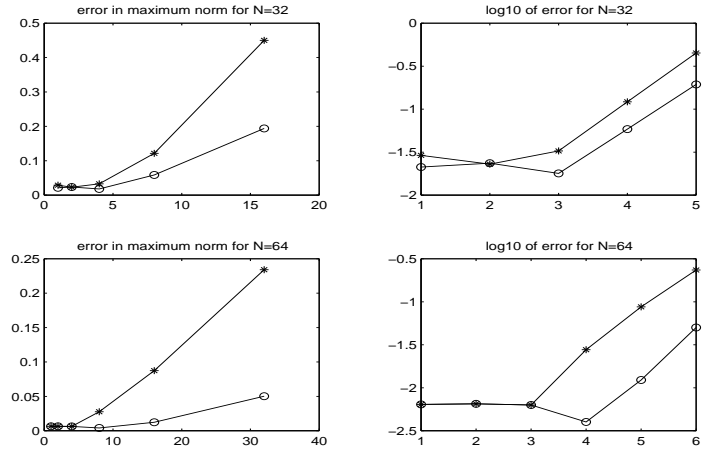,height=2.42in,width=3.7in}}
  \caption{Accuracy of the stabilized explicit scheme
   for the heat equation.
   The horizontal coordinate is $3\Delta t/h^2$.
   *~:~first-order shift; o~:~third-order shift.
  \label{GKR-fig1}}
  \end{center} 
  \end{figure}

Figure~\ref{GKR-fig1} shows some accuracy results
for Eq.~(\ref{GKR-scalar_reaction_diffusion}),
where
\begin{equation}
  u(x,t) = \cos(t) ( (x/\pi)^4 + \cos (3x)) , \quad
  x \in (0, \pi), \, t > 0 .
\end{equation}
\noindent
We observe a plateau for small time steps,
when the second-order spatial error dominates.
The second-order error in time becomes dominant 
as the time step increases.
The figure confirms the superior performance of
the third-order shift~(\ref{GKR-shift4})
over the first-order shift~(\ref{GKR-step1})
at large time steps.

\begin{figure}[htb]
  \begin{center}
  \mbox{\psfig{file=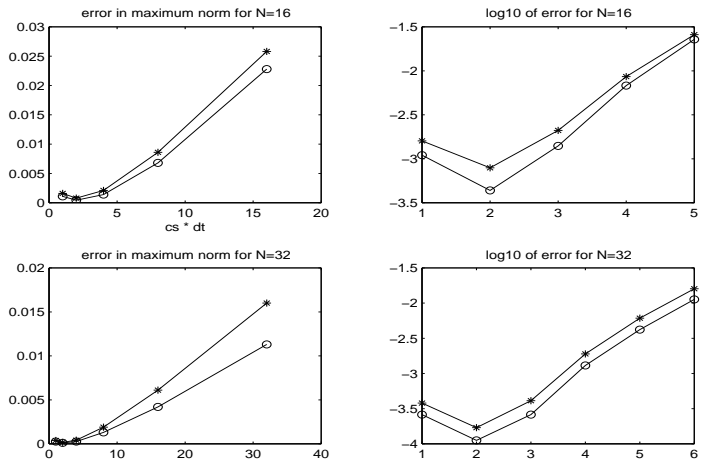,height=2.42in,width=3.7in}}
  \caption{Accuracy of the stabilized explicit scheme
   for the predator-prey system.
   The horizontal coordinate is $3 \Delta t/h^2$.
   *~:~first-order shift;
   o~:~third-order shift.
  \label{GKR-fig2}}
  \end{center} 
  \end{figure}

Although the algorithm is based only on linear
stability considerations,
it is still effective for systems of
nonlinear reaction-diffusion equations.
In Figure~\ref{GKR-fig2} we present some results for
a predator-prey system,
\begin{equation}
  \partial_t u = \partial_{xx} u + a u - b u v , \;
  \partial_t v = \partial_{xx} v - c u - d u v , \quad
  x \in (0,\pi), \, t > 0,
\end{equation}
with $a=1.2$, $b=1.0$, $c=0.1$, and $d=0.2$.
\noindent
At these parameter values,
the ODE system (reactions only) has a limit cycle.
However, when the boundary conditions are constant in time,
the solution of the reaction-diffusion system
goes to steady state.
To build a relevant test case for the algorithm,
we impose periodic excitations at both boundaries,
\[
  u (0 / \pi, t) = u_{0 / \pi} (1+\cos(t)), \;
  v (0 / \pi, t) = v_{0 / \pi} (1+\cos(t)), \quad t > 0 .
\]
Although the time step can still be limited
by nonlinear instabilities,
we never observed negative values of
the unknowns $u$ and $v$, which are
commonly associated with such instabilities.

\begin{figure}[htb]
  \begin{center}
  \mbox{\psfig{file=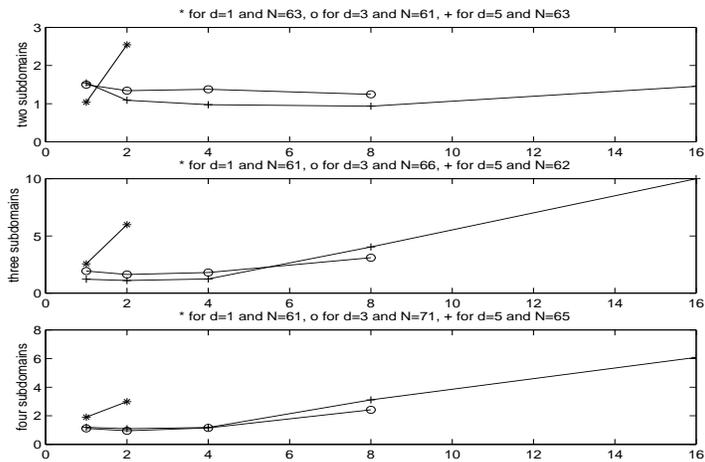,height=2.42in,width=3.7in}}
  \caption{Accuracy of the stabilized explicit scheme
   applied to the heat equation
   with overlapping subdomains.
  \label{GKR-fig3}}
  \end{center} 
  \end{figure}

The algorithm~(\ref{GKR-step1})--(\ref{GKR-step3})
extends in a straightforward way when one uses
a domain-decomposition scheme with overlapping subdomains.
One simply applies the same algorithm
at each time step to each subdomain separately.
However, the number of waves per subdomain
is of the order of the total number
of grid points, $N$, divided by the number of
subdomains, $N_d$, so the balance between
the order of accuracy of the filter---
$(N/(\kappa N_d))^{-3}$ for a first-order filter or
$(N/(\kappa N_d))^{-5}$ for a third-order filter---and
the second-order accuracy $N^{-2}$ of the spatial
discretization of the underlying
algorithm~(\ref{GKR-schema})
deteriorates as $N_d$ increases.
The maximum time step for which the scheme
remains stable may become even less
than when no domain decomposition is used.
Furthermore, the Gibbs phenomenon
tends to destabilize the algorithm.
This phenomenon is a consequence of the jump in the
derivatives at the endpoints of the subdomains
(second-order derivatives in the case of the
first-order shift~(\ref{GKR-step1}),
fourth-order derivatives in the case of the
third-order shift~(\ref{GKR-shift4})).
Since the Gibbs phenomenon arises at
the artificial interface and is damped
away from it, an increase of the overlap
generally produces a composite signal $u$
that has fewer oscillations than each of
the piecewise (overlapping) components.
One can therefore obtain good results
by adapting the size of the overlap.
The larger the overlap, the larger
the time step that can be taken; see Figure~\ref{GKR-fig3}.

\section{Two-Dimensional Domains}
We now consider a Dirichlet problem
in two dimensions,
\begin{eqnarray} \label{GKR-reaction_diffusion}
  & \partial_t u = \Delta u + f(u) , \quad (x,y) \in (0,\pi)^2, \, t>0 , & \\
  & u(x,0 / \pi)= g_{0 / \pi}(x), \; u(0 / \pi,y)= h_{0 / \pi}(y), \quad
  x,y \in (0,\pi) , &
\end{eqnarray}
where the functions $g$ satisfy
the compatibility conditions
$g_{0/ \pi}(0)=h_0(0 / \pi)$ and
$g_{0/ \pi}(\pi)=h_{\pi}(0 / \pi)$.
We consider a numerical scheme similar
to Eq.~(\ref{GKR-schema}), where
the diffusive term is approximated,
for example, by a five-point stencil.
The postprocessing algorithm is essentially the same,
except that we need an apropriate
low-frequency shift so we can apply
a filter to a smooth periodic function
in two space dimensions.
The shift is constructed in two steps.
In the first step, we render
the boundary condition in the $x$ direction
homogeneous,
\begin{eqnarray} \label{GKR-shift_x}
  &v(x,y)
  = u(x,y) - (\alpha_1 (y) + \alpha_2 (y) \cos(x)) , &\\
  & \alpha_1 (y) = \half (g_0 (y)-g_{\pi} (y)) , \,
    \alpha_2 (y) = \half (g_0 (y)+g_{\pi} (y)) . &
\end{eqnarray}
In the second step, we shift in the $y$ direction,
\begin{eqnarray} \label{GKR-shift_y}
  &w(x,y) = v(x,y) - (\beta_1(x) + \beta_2 (x) \cos(y)) , &\\
  &\beta_1 (x) = \half (v(x,0)-v(x,\pi)) , \,
   \beta_2 (x) = \half (v(x,0)+v(x,\pi)) . &
\end{eqnarray}
The final step is the reconstruction step,
\begin{equation} \label{GKR-step3_2D}
  u(x, y) = \sigma_N w(x, y)
  + \alpha_1 (y) + \alpha_2 (y) \cos(x)
  + \beta_1 (x) + \beta_2 (x) \cos(y) .
\end{equation}
To make sure that no high-frequency waves remain,
we filter the high-frequency components
from the boundary conditions $g$
with a procedure similar
to~(\ref{GKR-step1})--(\ref{GKR-step3}).

It is much more difficult to construct a high-order filter
similar to~(\ref{GKR-shift4}) in two dimensions,
because the second-order derivatives cannot be obtained
from the PDE, as in the one-dimensional case~(\ref{GKR-uxx}).
So far, we have used only the first-order
shifts~(\ref{GKR-shift_x}) and~(\ref{GKR-shift_y})
in our numerical experiments.
Nevertheless, the algorithm allows for
a significant increase of the time step.
We have also tested the domain-decomposition version
of the algorithm, using strip subdomains with
an adaptive overlap, with good results.

We note that the computation in each block
can be done in parallel and that the
Jacobi matrix does not depend on the spatial
variables.
The arithmetic complexity of the algorithm
is therefore relatively small.
Also, the algorithm is suitable
for multicluster architectures.
Each block can be assigned to a cluster,
and parallel fast sine transforms
can be used for the filtering process
inside each cluster.
The cost of communication between blocks
is minimal, since the scheme is similar
to the communication scheme
of the additive Schwarz algorithm.

\section{Conclusion} 
In this paper we have presented
a postprocessing algorithm
that stabilizes the time integration
of systems of reaction-diffusion equations
when the diffusion term is treated explicitly.
The algorithm is easy to code and
can be combined with domain-decomposition
methods that use regular grids in each subblock.
In future work, we will consider the performance
of its parallel implementation and its robustness
for large systems of reaction-diffusion equations
with stiff chemistry,
which arise in some air pollution models.

\bibliographystyle{unsrt}
\bibliography{GKR-dd13}
       
\end{document}